\documentclass{amsart}

\usepackage{amsmath, amsthm, amscd}
\usepackage{amssymb, latexsym}
\usepackage{amsfonts}

\newtheorem{theorem}[equation]{Theorem}

\newtheorem{proposition}[equation]{Proposition}
\newtheorem{lemma}[equation]{Lemma}

\newtheorem{corollary}[equation]{Corollary}

\newtheorem{Remark}[equation]{Remark}
\newtheorem{Remarks}[equation]{Remarks}

\newtheorem{Example}[equation]{Example}

\newcounter{com}

\newtheorem*{A'}{Corollary A$'$}

\newtheorem*{C}{Theorem C} 

\newtheorem*{C'}{Theorem C'}
\newtheorem*{C''}{Theorem B$''$ (Holt's Conjecture)}

\newtheorem*{CC}{Conjecture C}

\newcommand{\beql}[1]{\begin{equation}\label{#1}}
\newcommand{\eeq} {\end{equation}}

\font\Aaa=msam10

\def\qed{\hbox{~~\Aaa\char'003}}

\font\Bbb=msbm10

\def\W{\hbox{\Bbb W}}

\def\UU{\mathcal U}

\numberwithin{equation}{section}

\let\define=\def
\let\redefine=\def

\def\Char{{\rm char}}

\def\PG{{\rm PG}}
\define\a{{ \alpha }}
\redefine\b{{ \beta }}

\define\P{{ \bf P}}

\define\({ \left( }
\define\){ \right) }
\define\[{ \left[ }
\define\]{ \right] }
\define\<{ \langle }
\define\>{ \rangle }
\let\ljunk=\{
\let\rjunk=\}
\redefine\{{\left\ljunk}
\redefine\}{\right\rjunk}

\redefine\.{ \diamomd }
\redefine\+{\oplus}
\redefine\.{ \cdot}
\redefine\#{\sharp}

\redefine\.{ \cdot }

        \def\GF{{\rm GF}}
        \def\GL{{\rm GL}}

\def\P{{\bf P}}

        \font\Aaa=msam10

\begin{document}
\title[Veroneseans, power subspaces and independence]
{Veroneseans, power subspaces and independence}\thanks{This research was supported in part by
 NSF grant DMS~0753640.}

       \author{W. M. Kantor}
       \address{University of Oregon,
       Eugene, OR 97403}
       \email{kantor@uoregon.edu}

    \author{E. E. Shult}
       \address{Kansas State Univesity, Manhattan, KS  66502}         \email{ernest$\underline{~}$shult@yahoo.com}

 
\begin{abstract}  
 
Results are proved indicating that the Veronese map  $v_d$  often increases independence
of both  sets of points  and sets of subspaces.  
For example,  any $d+1$  Veronesean points of  degree $d$  are independent.  Similarly, the $d$th power map on the space of linear forms of a polynomial  algebra also  often increases independence
of  both sets of points and sets of subspaces.
These ideas produce $d+1$-independent families of subspaces
in a natural manner.  
 \end{abstract}
 
\maketitle

\vspace{-6pt}

\section {Introduction}
\label {Introduction}
 In this paper we will study independence questions involving 
 points or subspaces obtained from standard geometric or algebraic objects: Veronese maps and polynomial rings. 
  The proofs are elementary, but some of the results seem unexpected.
 
 We will always be considering   integers $d,n >1$.
For any field $K$, the 
(vector) {\em Veronese map}%
 \footnote{Traditionally this map is defined sending   
 the projective geometry $\P( K^n )$ on $K^n$ to $\P( K^{N} ).$
For our proofs it is preferable to deal with vectors, 
which still allows us to act on projective geometries.} 
 $v_d\colon\! K^n\to K^N$, 
 $ N=\binom{d+n-1}{d} , $ 
is defined in (\ref{Veronesean map}); 
\vspace{.5pt}
the $1$-subspaces in 
$v_d  ( K^n ) $  are the
 {\em Veronesean points of degree $d$}.  We will be concerned with 
 the  behavior of $v_d$ on  sets of  subspaces of $ K^n $: in general it increases independence.
 For example:

\begin{theorem}
\label{Veronesean point independence}
Any $d+1$ Veronesean  points of degree $d$  in $ K^{N}$
are independent  $($that is$,$ they span a $d+1-$space$)$.
\end{theorem}
The dimension $n$ of the initial space $K^n$
 does not play any role in this result or others in this
 paper.
Section~\ref{Veronesean points} contains a surprisingly elementary proof.  
These types of results are 
in the geometric framework appearing  in 
  \cite{Lu,HT,CLS,BC,KSS1,KSS2} rather than the more standard 
  Algebraic Geometry occurrences of the Veronese map
\cite[p.~23]{Har}, \cite[pp.~40-41]{Sh}.   
  
  More generally, we will 
 consider independence of sets   of subspaces of $K^n$. We call a 
  set $\mathcal U$ of at least $d+1$ such 
  subspaces  {\em$d+1-$independent}   if the subspace spanned by any $d+1$ members of $\mathcal U$ is their direct sum.  For example, $2-$independence means that any two members have intersection 0, which is a very familiar geometric situation.  
   With this terminology,
   Section~\ref{Families of subspaces} contains an elementary 
      proof of   the following generalization of the preceding theorem concerning sets $v_d({\mathcal U}):= \{\< v_d(U )\>\mid U\in {\mathcal U  }\}$ of subspaces $K^N$:
  
  \begin{theorem}
\label{Veronesean e+1 independence}
If $~\mathcal U$ is any $e+1-$independent set 
 of at least $de+1$ nonzero
 subspaces of $K^n,$ then 
 $v_d({\mathcal U})$ is a $de+1-$independent  set of subspaces of  $K^N$.
\end{theorem}

   Even when  $ \mathcal U$ is a spread we suspect 
%   \marginpar{IS THIS OK?}
that  the resulting $d+1$-independent family $ v_d(\UU)$ is not maximal.  It is the  natural way of obtaining $ v_d(\UU)$   that seems  more interesting   than the  possible maximality.   Note that the dimensions of the subspaces in the preceding theorem are allowed to vary arbitrarily.

For any   finite set  $\{ x_1,\ldots , x_n\}$
of indeterminates and any integer
 $m\ge1$,  we will consider the space $A_m$ consisting of all   homogeneous polynomials of degree $m$  in the polynomial algebra $A=K[x_1,\ldots , x_n]$.  A {\em powerpoint} is a 1-space $\<f^d\>$
in $  A_d$, where $0\ne f\in A_1$.  Powerpoints in $A_d$ are 
closely related to  Veronesean points in suitable characteristics (cf. Theorem~\ref{th:pover}).  As in the case of the Veronese map, we are interested in the behavior of $d$-fold powers on  sets of subspaces of 
$ A_1$.   The case of powerpoints follows easily from 
 Theorem~\ref{Veronesean point independence}  (cf.  Section~\ref{Powerpoints and Veronesean points}):

\begin{theorem}
\label{power point independence}
For any field whose characteristic is $0$ or $\,>d, $  any $d+1$ powerpoints of $A_d$ 
are independent.
\end{theorem}
 
 More specialized results are possible for small positive characteristics
(cf. Theorems~\ref{powerpoint independence  small char} and \ref{powerpoint independence char 2}).

Let $\<T ^d\>$ denote the subspace spanned by all products of $d$ members of a subset $T$  of  $ A$.
  Section~\ref{Proof of Theorem independence} again  concerns increasing independence of subspaces, once again assuming a restriction on the characteristic:

\begin{theorem} 
\label{independence}
Assume that $1\le r\le d$ and  $d!/(d-r)!\ne0$ in $K$. 
If  $~\mathcal T$ is any $e+1-$independent set 
  of at least $re+1$ nonzero
  subspaces of $A_1,$ then $\{ \<T ^d\>\mid T\in {\mathcal T }\}$ is an $re+1-$independent set in $A_d$.
\end{theorem}  
Theorem~\ref{Veronesean e+1 independence}  can
 be used to prove this when $r=d$, while Theorem~\ref{power point independence} is a special case, although the proofs use  very different tools.
In Section~\ref{3-Independence of  power subspaces} we prove a somewhat weaker-looking
variation on the preceding theorem.

Generalized dual arcs and other configurations are constructed in Section~\ref{Generalized dual arcs and polynomial rings}  using very elementary properties of the polynomial algebra $A$.  
One of these configurations is  another infinite family of $3$-independent subspaces.

Note:  We will always use vector  space dimension.

\section{The Veronese map}
\label{The Veronesean map}

In this section we will prove Theorems~\ref{Veronesean point independence} and \ref{Veronesean e+1 independence}.  In passing we use the polynomial ring to reprove a standard result on Veronesean action.

\subsection{Background concerning the Veronese map}
\label{Background concerning the Veronese map}
Consider two integers $d, n>1$,
 together with   $N=\binom{d+n-1}{d} $.
  For a field $K$ of arbitrary 
  \vspace{2pt}
  characteristic and size, we will use the $K$-space 
    $V=K^n$ of all $n$-tuples $t=(t_1,\dots,t_n)=(t_i),$ $ t_i\in K$, and 
     the $K$-space $W=K^N$ of all $N$-tuples 
 $(y_\a)$   for a fixed but arbitrary ordering of all sequences $\a=(a_1,\dots,a_n)$ of integers $a_k \ge0$ satisfying $\sum_ka_k=d$.
 Corresponding to $\a$ there is a monomial function
  $t\mapsto t^\a:=t_1^{a_1}\cdots t_n^{a_n}$ of degree $d$ in the coordinates $t_i$.  (See Section~\ref{Background}   for more discussion of 
  this setting.)

The (vector) {\em Veronese map}  $v_d\colon \!V\to W$
is defined by
\begin{equation}
\label{Veronesean map}
v_d \big((t_i) \big)  := (t ^\a).
\end{equation}
This   induces the classical Veronese  map  $\P(V)\to \P(W)$ on projective spaces \cite[p.~23]{Har}, \cite[pp.~40-41]{Sh}.
Some of its geometric aspects have been studied outside Algebraic Geometry in 
  \cite{Lu,HT,CLS,BC,KSS1,KSS2}.

 There is a natural map from homogeneous polynomial 
 functions $g$ in
$t_1,\dots,t_n$ of degree $d$ to linear functionals $W\to K$. Namely, if $g(t_1, \dots , t_n) = \sum_\a  a_\a t^\a$, where
$a_\a\in K$ using all $\a$ as before, then the corresponding
linear functional  $\tilde g\colon\! W\to K $ is given by 
 $\tilde g((y_\a)) := \sum_\a  a_\a y_\a $. If the field is tiny then it is possible that two monomial functions $t^\a$ coincide, so that this correspondence is not bijective, in fact ``sending'' $g$ to
  $\tilde g$  is not actually a function!  However, what matters here  is that this recipe produces a linear functional, and  that every 
 linear functional on $W$ arises this way.

 Clearly, $\tilde g$ is a linear functional on $W$ such that
\begin{equation}
\label{g tilde}
 \mbox{$\tilde g\big ( v_d(t)  \big) = g(t)$   for all $t= (t_i)\in  V$.}
\end{equation}

 \subsection{Veronesean points}
 \label{Veronesean points}

The following elementary observation implies Theorem \ref{Veronesean point independence} (see  Lemma~\ref{almost independence} for a much stronger  version):

\begin{lemma}
\label{Veronesean subspaces}
 If $z$ is a point of $V$ not in each of $d$
  subspaces $U_1, \dots,U_d$ of $\,V ,$
then $v_d(z)$ is not in 
$\<v_d (U_1), \dots  ,v_d (U_d) \>$.
\end{lemma}

\proof 
 For each $j$ let $f_j$ be a linear function   $V\to K$ that vanishes on 
 $U_j$ but not on $z$.
Then $g := \prod_j f_j$ is a homogeneous polynomial function of degree $d$ that vanishes
on all $U_j$ but not on $z$. By (\ref{g tilde}), the corresponding linear function $\tilde g$ on $W$ vanishes
on $\<v_d (U_1), \dots  ,v_d (U_d) \>$ but not on $v_d(z)$, as required.
 \qed 
 
%%% \newpage
 
\medskip
  
See  \cite{BC} and \cite[Theorem~2.10]{CLS} for results similar to  Theorem~\ref{Veronesean point independence}.

\Remark \rm
 \label{code}
If $q\ge d$ then the  rational normal curve $v_d(\P(K^2))$ spans $K^N=K^{d+1}$ 
\cite[p.~229]{Hi}.  It follows that $v_d(K^2)$ does not contain $d+2$ independent points:  Theorem~\ref{Veronesean point independence}
is best possible.

Theorem~\ref{Veronesean point independence} can be viewed as a statement about the code $C$ having a check matrix whose columns 
consist of one nonzero vector in each Veronesean point:   
 $C$ has minimum weight $>d+1$.  By the preceding paragraph, the minimum weight  is $d+2$ if $q$ is not too small, with codewords of weight $d+2$ arising from $d+2$ points in a $2-$space in $K^n$;
 and similarly, the next smallest weight is $2d+2$, occurring from  
 $ d+1$ points on each of two $2-$spaces  in a $3-$space. 
 
 We have not been able to find any reference to this code in the literature.  It is probably worth studying, at least  from a geometric perspective.
\Remark \rm
 \label{iterate}
The notation $v_d$ is ambiguous, since it omits the original 
dimension $n$.  With this in mind, these maps can be composed.
It is easy to use monomials to check that $v_e(v_d(\P(K^n)) )$
 is just $v_{ed}(\P(K^n))$ 
    on a subspace of  the underlying $\binom{e+N-1}{e} $-dimensional space (where $N$ is as before).
 For example, for any set $X$ of points of a projective space, 
 $v_e(\P(X))$ is $e+1-$independent by 
Theorem~\ref{Veronesean point independence}; but if $X=v_d(\P(Y))$  for $Y\subseteq A_1$  then
 $v_e(\P(X))$ is   $de+1-$independent.
 
 In particular, if $C$ is a conic in $K^3$ then  
$v_2(C)
=v_2(v_2(\P(K^2)))$ is 5-independent: it is a rational normal curve
in a 5-space.
Similarly, it is  natural to ask for the independence properties of $v_d$-images of geometrically natural sets of points.
For example, by Theorem~\ref{Veronesean e+1 independence}, 
$v_2(\mbox{hyperoval})$
and $v_2(\mbox{ovoid})$ are 5-independent (and 5 is best 
possible). 
    
 \subsection{Families of subspaces}
 \label{Families of subspaces}
 The following  special case of Theorem~\ref{Veronesean e+1 independence} contains Theorem~\ref{Veronesean point independence}:
   
\begin{theorem}
\label{Veronesean independence}
 If $~\mathcal U$ is any set 
 of at least $d+1$  nonzero
 subspaces of $K^n$   pairwise intersecting in $0,$ then 
 $\< v_d({\mathcal U  })\>$ is a $d+1-$independent  set 
 of subspaces of  $K^N$. 
\end{theorem}
 
This is an immediate consequence of the following

\begin{lemma}
\label{almost independence}
 If $U_0$ is a subspace of $V$ intersecting each of $d$
  subspaces $U_1, \dots,U_d$ of $\,V  $ only in $0,$
then $\<v_d(U_0)\>\cap \big\<v_d (U_1), \dots ,v_d(U_d) \>=0.$
\end{lemma}

\proof We will construct a linear map $L$  on $ W$ 
whose kernel contains $v_d (U_j)$, 
$1\le j\le d$,  and meets $\<v_d(U_0)\>$ only in 0.
By Corollary~\ref{cor:rho}, one  may arbitrarily change a basis of $V$ while leaving the set $v_d(V)\subseteq W$ invariant.  Thus we may assume that
$U_0=\{(t_1 , \dots,t_m,0,\dots, 0)\mid t_i\in K\}$;
we will view $U_0$ as $K^m$.  
Let $W_0\cong K^{N_0}$,
$N_0=\binom{d+m-1}{d} $,
 be the set of vectors in $W$ having a nonzero coordinate for some member of $v_d(U_0)$.  Then 
$v_d\colon \! U_0\to W_0$ \emph{can be viewed 
as the Veronese map on $U_0$.} (If  $K$ is small then $W_0$ might not be the span of $v_d(U_0)$, which  adds a minor complication
to  our argument.)

Throughout  this proof, $\b$ will
range over all sequences $(b_1,\dots,b_m,0,\dots,0)$ of $n$  integers $b_k \ge0$ satisfying $\sum_kb_k=d$.

Let ${\mathbb H}$ denote the $K$-space of all homogeneous polynomial functions $g$ on $V$ 
of degree $d$ such that $g(U_j)=0$ for $1\le j\le d$. We will construct many elements of  ${\mathbb H}$. 
First    note that \emph{every monomial function
$t^\b$ on $U_0$ of degree $d$ is the restriction of some member of} ${\mathbb H}$. 
For, write $t^\b = t_{\sigma(1)}\cdots t_{\sigma(d)}$ for a function
$\sigma\colon \!\{ 1,\dots,d\}\to  \{ 1,\dots,m\}$.  
If  $1\le i\le d$   
 let $\lambda_i$ denote any linear functional on $V$ such that 
$\lambda_i \big((t_1 , \dots,t_m,0,\dots, 0) \big)=t_{\sigma(i)}$ and 
$\lambda_i(U_i)=0$ (recall that $U_0\cap U_i=0$).  Then   $ \prod_{i=1} ^d\lambda_{i}
\big((t_1 , \dots,t_m,0,\dots, 0) \big) 
= \prod_{i=1} ^d t_{\sigma(i)} =t^\b $ and 
$ \prod_{i=1} ^d\lambda_ i (U_j)=0$  for $1\le j\le d$. Thus, 
$\prod_{i=1} ^d\lambda_i $ behaves as required.

It follows that every homogeneous polynomial function on $U_0$ of degree $d$ is  the restriction of some member of ${\mathbb H}.$

Let $\W$ denote the set of all linear functionals on $W$ that vanish on
$v_d(U_j)$ for $1\le j\le d$.  
Set $W_{\bullet}:=\<v_d(U_0  )\>$.  
The crucial step of the proof of the lemma is that
\begin{equation}
\label{crucial}
\mbox{\em Every 
linear functional $\mu_{\bullet}$ on $ W_{\bullet}$   is  the restriction 
of some  $\mu \in \W$.}
\end{equation}
For, arbitrarily  extend  $\mu_{\bullet}$  to a linear functional $\mu_0$ on $W_0$ (this is irrelevant if $W_0=W_{\bullet}$).
As noted in Section~\ref{Background concerning the Veronese map}, there is a 
homogeneous polynomial function 
$g_0$ on $U_0$ of degree $d$ such that $\mu_0=\tilde g_0$.
We have seen that $g_0$ is the restriction to $U_0$ of some 
$g\in {\mathbb H}$.
Consequently, it suffices to show that $\mu:=\tilde g$ coincides with 
$\mu_{\bullet}$ on $W_{\bullet}$.  If $t\in U_0$ then we can apply (\ref{g tilde}) using both $V$ and~$U_0$:
$$
\tilde g \big(v_d(t) \big)= g(t)=g_0(t)=\tilde g_0 \big(v_d(t) \big) =\mu_0 \big(v_d(t) \big)=\mu_{\bullet} \big(v_d(t) \big).
$$
Since $\tilde g$ and $\mu_{\bullet}$ are linear on $W_{\bullet}=\< v_d(U_0) \>, $
it follows that  $\tilde g = \mu_{\bullet}$ on $W_{\bullet}$,
which proves (\ref{crucial}).

%\bigskip \bigskip
 
 Set   $N_{\bullet}:=\dim W_{\bullet}$.  
 By (\ref{crucial}), $\W$  has a subset  $\{\mu_i\mid 1\le i\le N_{\bullet}\}$  whose restrictions to $W_ {\bullet} $ form
 a basis of the dual space $W_ {\bullet} ^{*}$.
 Then $\mu_i \big(v_d(U_j) \big)=0$ for $1\le i\le N_{\bullet},$  $1\le j\le d$,
 by the definition of $\W$.

Define $L\colon W\to K^{N_{\bullet}}$ by 
$L \big((y_\a) \big)= \big(\mu_i \big((y_a) \big) \big)$.  
Then $L$ is linear, and
$L \big( v_d(U_j) \big)= \big(\mu_i \big(v_d(U_j) \big) \big)=0$  for $1\le j\le d$. 
Since $\{\mu_i\mid 1\le i\le N_{\bullet}\}$ restricts to a basis of $W_{\bullet}^{*} $,
 % we have  
$$\hspace{24pt}\<v_d(U_0)\>\cap \big\<v_d (U_1), \dots ,v_d(U_d) \>
\le W_{\bullet}\cap \ker L = W_{\bullet}\cap \bigcap_i\ker \mu_i=0.
 \hspace{24pt}
\qed$$
  
  %\medskip
  
\Remark \rm
The most familiar examples of $2$-independent families are spreads.
It would be interesting to know for which $r$ the set in 
Theorem~\ref{Veronesean independence}
is $r-$independent 
when $\Sigma$ is a Desarguesian spread of $k$-spaces of a $2k$-space. 
\medskip

  The preceding lemma also yields {\rm Theorem~\ref{Veronesean e+1 independence}}:
  \medskip

{\noindent\em Proof of\: {\rm Theorem~\ref{Veronesean e+1 independence}}.}
Consider distinct  $U_0,\dots ,U_{de}\in {\mathcal U  }$,
and suppose that $\sum_{i=0}^{de}y_i=0$ for some $y_i\in  \<v_d(U_i) \>$.
By symmetry, it suffices to show that $y_0=0$.

Let $\Pi$ be any partition of $\{1,\dots,de\}$ 
into $d$ subsets $\pi$ of size $e$.
For $\pi\in \Pi$
let $U_\pi:=\<U_i\mid i\in \pi\>$.  Then 
$U_0\cap U_\pi=0$ since $\{U_0, U_i\mid i\in \pi\}$ is $e+1$-independent.  By the preceding lemma,
$$
\begin{array}{llll}
-y_0=\sum_1^{de}y_i\hspace{-6pt}&\in & \hspace{-8pt}\displaystyle
\<v_d(U_0)\>
\cap 
\sum_{\pi\in\Pi}
\sum_{i\in \pi}\<v_d(U_i) \> 
\vspace{4pt}
\\
&\le& \hspace{-6pt}
\<v_d(U_0)\>\cap 
\< v_d(U_ \pi)  \mid \pi\in\Pi \>
=0. \hspace{24pt}\qed
\end{array}
$$

% \big\medskip

\Remarks  \rm 1. We used the rather weak inclusion 
$\<v_d(A), v_d(B )\>\le \<v_d(\<A,B \>)\> $ for subspaces $A,B$ of $K^n$:  in general the right side is far larger than the left.
\smallskip

2. The proof shows that we did not need independence for all 
$e+1$-subsets of $\mathcal U$. For each $de+1$-subset $\mathcal U'$
of $\mathcal U$ we only needed a family $\W$ of independent $e+1$-subsets of $\mathcal U'$  such that the complement of each member of $\mathcal U'$ is partitioned by some of the members of $\W$.

The minimal version of this is as follows: each $de+1$-subset $\mathcal U'$
of $\mathcal U$ is equipped with the structure of a 2-design with $ v=de+1, k=e+1, \lambda=1,$ such that each block is $e+1-$independent.
In this situation,  ``almost all''   triples  from  $\UU$ need \emph{not} be independent and yet the proof shows that $v_d(\UU)$ 
nevertheless must be $de+1-$independent.

\subsection{Veronesean action}
\label{action} 
\label{Background}
This section develops two algebraic results that play a small role in the proofs in this paper. 
 One is that   linear transformations of the space of homogeneous polynomials of degree one 
induce  endomorphisms of degree zero of the polynomial algebra $K[X]$   (see Remark~\ref{extends to endomorphism}).  The other is the oft-quoted result that there is an action of $\GL(K^n)$ on $K^N$ that 
stabilizes the set of Veronesean vectors, inducing   an action permutation-equivalent to its action on $K^n$ (used in Lemma \ref{almost independence}), which we prove  
using  polynomial rings and their morphisms. 
\subsubsection{Symmetric algebras and polynomial rings}

Let $V$ be an arbitrary vector space over   
 $K$ of dimension $n$.  The symmetric algebra $S(V)$ is the $K$-algebra of symmetric tensors -- that is, the free commutative $K$-algebra generated  by vector space $V$.  It is a graded algebra 
$$S(V) = K\oplus V\oplus S_2(V)\oplus \cdots $$
where $S_d(V)$ is the vector space spanned by the $d$-fold symmetric tensors.  If $X = \{ x_1, x_2,\ldots , x_n\}$ is any basis of   $V$
then  $S(V)$ is isomorphic to the polynomial ring 
\begin{equation}
\label{eq:graded}
A = K[X] = K\oplus A_1\oplus \cdots \oplus A_d\oplus\cdots 
\end{equation}
where $A_d$ is the vector space of homogeneous polynomials of degree $d$.  Thus selecting the basis 
$X$ of $A_1$ 
produces a basis $\{ x^{\alpha}\}$ of $A_d$ consisting of
the  monomials $x^{\alpha}:= x_1^{a_1}\cdots x_n^{a_n}$, where $\alpha = (a_1,\ldots , a_n)$ is a sequence of non-negative integers for which $d = \sum a_i$.

\subsubsection{The substitution-transformation $\rho_d$}
Let $f:A_1\rightarrow W$ be any linear transformation, where $W$ is
 a $K$-vector space.  Then $f$ extends to a $K$-algebra homomorphism $\bar{f}: A[X]\rightarrow S(W)$ of graded algebras, by mapping any polynomial $p(x_1, \ldots , x_n)$ to $p(f(x_1), f(x_2), \ldots , f(x_n))$, a ``polynomial'' in the algebra $S(W)$.  By restriction of $\bar{f}$, we set 
$$ \rho_d(f): = \bar{f}|_{A_d}: A_d\rightarrow S(W)_d.$$
Its value at any monomial $x^{\alpha} = \prod x_i^{a_i}$ is $\prod f(x_i)^{a_i}$ in $S(V)_d$.
  Thus $\rho_d(f)$ is simply the linear morphism on $A_d$, which results from {\em substituting} $f(x_i)$ for $x_i$.  
\Remark\label{extends to endomorphism}\rm Note that when
$W\leq A_1$, $f$ has been extended to an {\em endomorphism} of the algebra $K[X]$.

Now what happens when we apply $\rho_d$ to a functional $\lambda : A_1\rightarrow K$?  Since $\bar{\rho}(\lambda )$ is defined by substitution of each $x_i$ by the scalar $\lambda (x_i)$ in each polynomial of $K[X]$, it induces a functional $\rho_d(\lambda): A_d\rightarrow K$ of $A_d$.
\begin{lemma}\label{lem:rho}
Some properties of $\rho_d$:
\begin{enumerate}
  \item $\rho_d$ transforms any linear transformation  $A_1\rightarrow A_1$ to a linear transformation of $A_d$ into itself.  If $T$ is the identity  transformation of $A_1$, then $\rho_d(T)$ is the identity transformation of   $A_d$.
   \item If $T$ is in the group $GL(A_1)$, then $\rho_d(T)$ is an invertible transformation of $A_d$.
   \item If $\lambda:A_1\rightarrow K$ is a functional of $A_1$, then $\rho_d(\lambda)$ is a functional of $A_d$.
  \item If $S: A_1\rightarrow A_1$ is a linear transformation and if $T: A_1\rightarrow W$, where $W$ is either the $K$-vector space $A_1$ or the $K$-algebra $K$ itself,  then 
  \begin{equation} \label{eq:comp}
  \rho_d(T\circ S) = \rho_d(T)\circ \rho_d(S),
  \end{equation}
  and is also  $K$-linear.
  \item  Suppose $R$, $S$  are linear transformations $A_1\rightarrow A_1$ while $T: A_1\rightarrow W$ 
  is also $K$-linear, where $W$ is as in (4).  Then 
  \begin{equation}\label{eq:rho-comp2}
  \rho_d(T\circ S\circ R) = \rho_d(T)\circ\rho_d(S)\circ\rho_d(R)
  \end{equation}
 
\end{enumerate}
\end{lemma}

\proof
The first part of (1) follows from the fact that $\rho_d(T)$ is defined  by substituting $T(x_i)$ for $x_i$ in any homogeneous polynomial of degree $d$.  If $T$ is the identity map on $A_1$, then substitution of $x_i$ for $x_i$, does not change anything --- that is,  $\rho_d(T)$ is the identity transformation of $A_d$.

Statement (3) was explained in the paragraph preceding the lemma.

Statement (4) is also a consequence of $\rho_d(T)$ being defined by ``substitution".  Since $S$ is a linear transformation of $A_1$ into itself, we may utilize the basis $X$ to write 
$$S(x_i)= \sum_{j=1}^n c_{ij}x_j, \mbox{ where } c_{ij}\in K, i\in [1,n].$$
Then $\rho_d(S)$ takes monomial $x^{\alpha} = \prod x_i^{a_i}$ to $\prod (\sum_j c_{ij}x_j)^{a_i}$.  Since $\rho_d(T)$ takes any monomial $\prod x_j^{b_j}$ of degree $d$ to $\prod T(x_j)^{b_j}$, we see that
\begin{equation}\label{eq:rho-comp}
 \rho_d(T)\circ \rho_d(S): x^{\alpha} = \prod x_i^{a_i} \mapsto \prod_i(\sum c_{ij}T(x_j))^{a_i}.
 \end{equation}
But since $T\circ S$ takes $x_i$ to $\sum_j c_{ij} T(x_j)$ we see that  it also takes the monomial $x^{\alpha}$ to the right side of equation (\ref{eq:rho-comp}).  Thus we have
$$ \rho_d(T\circ S) = \rho_d(T)\circ  \rho_d(S),$$
establishing statement (4).

Remembering that $W$ is permitted to  be $A_1$ in statement (4), statement (5) follows from applying equation (\ref{eq:comp}) several times. 

For statement (2) suppose $T$ is invertible, so there exists a
 $T^{-1}:A_1\rightarrow A_1$ such that $T\circ T^{-1} = id_1$, the identity transformation of $A_1$.  Now by (1) and (2.14), the identity transformation $id_d$ of $A_d$ can be written as 
$$ id_d = \rho_d(id_1) = \rho_d(T\circ T^{-1}) = \rho_d(T)\circ \rho_d(T^{-1}),$$
proving that $\rho_d(T)$ is invertible.
~\qedsymbol

\subsubsection{Veronesean functionals}
Suppose $\lambda \in A_1^*$ is the functional $A_1\rightarrow K$ that takes the
basis element $x_i$ to the scalar $t_i$.  Then $\rho_d(\lambda)$ is the functional of $A_d$ that takes the basis element 
$x^{\alpha} $ to $t^{\alpha}\in K$.  We call a functional of this type (that is, one that maps $x^{\alpha}$ to $t^{\alpha}$ where $t = (t_1,\ldots , t_n)$) a {\em Veronesean functional} of $A_d$.  These are very special elements of $A_d^*$.  

\begin{theorem}\label{th:rho}  The group $\rho_d(GL(A_1))$  induces  an action on $A_d^*$ that stabilizes  the set of non-zero Veronesean functionals in $A_d^*$ and induces on this set an action that is permutation-equivalent to the action of $GL(A_1)$ on the non-zero vectors of $A_1^*$ .    Explicitly, if $T\in GL(A_1)$ acts on $A_1^*$ by sending the functional $\lambda$ to $\lambda\circ T$, then $\rho_d(T)$ acts on $\rho_d(\lambda)$, the corresponding Veronesean functional, by sending it to $\rho_d(\lambda)\circ \rho_d(T) = \rho_d(\lambda\circ T)$, another Veronesean functional.
\end{theorem}

\proof If  $\lambda \in A_1^*$ and $S$ and $T$ are elements of $GL(A_1)$, then by Lemma \ref{lem:rho}
\begin{equation}\label{eq:action}
\rho_d(\lambda\circ S\circ T) = \rho_d(\lambda)\circ \rho_d(S\circ T) = \rho_d(\lambda)\circ \rho_d(S)\circ\rho_d(T),
\end{equation}
for any $\lambda\in A_1^*$.

 By (\ref{eq:action}), we have a right action of $\rho_d(GL(A_1))$ on the set of Veronesean functionals.  Since these functionals are in one-to-one correspondence with the elements of $A_1^*$, the equation
 \begin{equation*}
 \rho_d(\lambda\circ T) = \rho_d(\lambda)\circ\rho_d(T)
 \end{equation*}
exhibits the permutation-equivalence of the action of $GL(A_1)$ on $A_1^*$ and the action of its isomorphic copy $\rho_d(GL(A_1))$ on the Veronesean functionals of $A_d^*$. ~\qedsymbol 

\subsubsection{The Veronesean action}
\label{The Veronesean action}

\begin{corollary} \label{cor:rho} There is an action of $GL(A_1)$ on the non-zero Veronesean vectors of $K^N$ that is permutation equivalent to its action on the non-zero vectors of $A_1^*$, or, equivalently, its action as $GL(K^n)$ on $K^n$.
\end{corollary}  

 \proof As before, $\alpha = (a_1,\ldots , a_n)$ is a sequence of non-negative integers summing to  $d$, so that $x^{\alpha} =\prod x_i^{a_i}$  is a monomial of degree $d$.  We define the scalar  $t^{\alpha} = \prod t_i^{a_i}$ whenever $t = (t_,\ldots , t_n)\in K^n$.  If $N$ is the number of monomials of degree $d$ in $n$ indeterminates, then the classical {\em  Veronesean vectors} are the $N$-tuples of the form $(t^{\alpha})$.  
Define the action of $GL(A_1)$ on $A_d^*$, by $f^T:= f\circ \rho_d(T)$, for every functional $f\in A_d^*$ and $T\in GL(A_1)$.  Equation (\ref{eq:comp}) shows that this meets the definition of a group action.  By 
 Theorem \ref{th:rho}, this action stabilizes the set of Veronesean functionals.    

The vector space isomorphism $\tau : K^N\rightarrow A_d^*$ which maps $(y_{\alpha})$ to the functional of $A_d$ whose value on $x^{\alpha}$ is $y_{\alpha}$, bijectively maps the 
set $v_d(K^n)$ of Veronesean vectors of $K^N$ to the set of Veronesean functionals of $A_d^*$.  Conjugation by $\tau$ then transports this action of $GL(A_1)$ on non-zero Veronesean functionals  described in the previous paragraph, to an equivalent action on the non-zero 
 Veronesean vectors.  

Similarly, let $\mu:  K^n\rightarrow  A_1^*$ be the vector space isomorphism which maps an $n$-tuple $(t_i)$ to the functional on $A_1^*$ whose value at
 $x_i$ is $t_i$.  Then conjugation by $\mu^{-1}$ transports the action of  $GL(A_1)$ on $A_1$ to an action as the full linear group on $K^n$.  
One can express this in terms of the (vector) Veronesean mapping introduced in (\ref{Veronesean map}).   Thus, setting $t = (t_i)\in K^n$, 
$$ \tau (v_d(t)) = \rho_d(\mu(t)).$$
Then for any $S\in GL(A_1)$, we have  
$$ v_d(t)^S:=  \tau^{-1}\rho_d(\mu(t))\circ\rho_d(g)\circ \tau = v_d(\mu^{-1}\circ g\circ \mu) := v_d((t^S)).$$
Equality 
of the  extremal members of this equation justifies the last remark of the Corollary.
~\qedsymbol 

\medskip
 See \cite[(2.3)]{Herzer} and \cite[Theorem~2.10]{CLS} for   other approaches to this corollary.

    \section{Powerpoints}
      For the rest of this paper, $ x_1,\ldots , x_n$ will denote indeterminates over  $K$, and $A:=   K[x_1,\ldots , x_n]$ is the  graded algebra  (\ref{eq:graded}),
so  that
$ A_i  A_j 
\subseteq A_{i+j} $ for all non-negative integers $i, j$.  
(If $P$ and $Q$ are sets of polynomials then {\em$P Q$ will denote the set of all products $pq, $ $ p\in P,$ $ q\in Q$}.  In general, it is not a subspace even if $P$ and $Q$ are.) 
If we replace $\{x_1,\ldots , x_n\}$ by any other basis of $A_1$ then we still obtain the same subspaces $A_d$ (cf. Section~\ref{Background}).

In this and the next two sections we will be concerned with {\em powers $U^d$ of subspaces  $U $  of $A_1$.} 
For now we will consider the set $P_d(A_1)$ of {\em powerpoints} $U^d$: the case in  which $U $ has dimension 1,
in which case  so does  $U^d$.
  
 \subsection{Powerpoints and Veronesean points}
  \label{Powerpoints and Veronesean points}
It is elementary and standard 
 that these two types of  points are closely related for suitable characteristics:

\begin{theorem}
\label{th:pover}  
If $  \Char K>d$ or   $\Char K =0,$    then   there is a linear isomorphism 
$\sigma\colon\! A_d\to K^{N} , N= \binom{d+n-1}{d}, $ 
such that
\begin{itemize}
\item [\rm(a)]
$\sigma$ sends the set 
 of powerpoints  in $A_d$  to the set   $v_{d}(\P (K^n))$ 
  of Veronesean points in $K^N,$ and 
\item [\rm(b)]
$\sigma \big (\big[\eta \big((t_i) \big) \big]^d \big)=(t^\a)$  if $\eta\colon \!K^n\to A_1$ 
sends $(t_i)\mapsto \sum_i t_ix_i$.
\end{itemize}

\end{theorem}
Here $(t^\a)$ was defined  in the preceding section.
\proof 
 By the Multinomial Theorem,
each powerpoint is spanned by a polynomial of the form
$$ (t_1x_1 + \cdots + t_n x_n)^d = \sum_{\a}  c(\a)t^\a x^\a $$
with $t_i\in K$  and multinomial coefficients $c(\a)$.   
All $c(\a)$ are nonzero in view of the assumed characteristic.  Hence, the map  $\sigma$ defined by 
$ \sigma\colon \! \sum_{\a}  c(\a)k_\a x^\a \mapsto (k_\a),$  $ k_\a\in K, $  behaves
as required.  \qed
\bigskip

\noindent{\em Proof of\/}  Theorem~\ref{power point independence}. 
   Theorem \ref{th:pover} shows that
    linear independence of powerpoints corresponds to linear independence of Veronesean points.   Now use  Theorem~\ref{Veronesean point independence}.~\qed
%\smallskip 

\subsection{Small characteristic}

We   now use Remark~\ref{extends to endomorphism}
to prove   additional independence results in small characteristics, a situation excluded in Theorem~\ref{power point independence}:

  \begin {theorem} 
   \label{powerpoint independence  small char}
   Assume that $r$ is such that
   $|K| >(r+1)^2/2$
 and    $\binom{d}{i}\ne0$ in $K$ for $0\le i\le r$.
Then any $r+1$ powerpoints  of $A_d$
are independent.
   \end {theorem}

\proof
If $n=2$ then all powerpoints are spanned either by  $x_1^d$ or
$(x_2+t x_1) ^d = 
\sum_0^{r}\binom{d}{i} t^ix_2^{d-i}x_1^i
+ \sum_{r+1}^{d}\binom{d}{i} t^ix_2^{d-i}x_1^i$ for some $t\in K$.
Since $\binom{d}{i}\ne0$ for $0\le i\le r$,
it suffices to note that the Vandermonde determinant $\det(t_j^i)_0^{r}\ne0$ for any $r+1$ different elements  $t_j\in K$.

If $n>2$, assume that the result holds for $n-1$ indeterminates $x_i$. 
Consider  $r+1$ distinct powerpoints  
 $\<f_1 ^d\>, \dots , \<f_{r+1} ^d\>$  and  a linear dependence relation $\sum_1^{r+1}k_if_i^d=0$,  $k_i\in K$.  
 Apply  the endomorphism $T$ of $A$ fixing 
 $x_1,\dots,x_{n-1}$ and sending  $x_n$ to an arbitrarily chosen  linear combination $f $ of 
 $x_1,\dots,x_{n-1}$
 (cf. Remark~\ref{extends to endomorphism}).  
 This produces an identity $\sum_1^{r+1}k_iT(f_i)^d=0$ in the ring
  $T(A) = K[x_1,\dots,x_{n-1}]$.  If the powerpoints 
$\<T(f_i)^d\>, 1\le i\le r+1,$ are distinct then 
  induction implies that all $k_i$ are~0.

If $\<T(f_i)^d\>=\<T(f_j)^d\>$ for some $i\ne j$ then 
$T\<f_i\>=T\<f_j\>$, so that $\<f_i\> $ and $ \<f_j\>$
are congruent modulo $\ker T=\<x_{n}-f\>$.  Therefore, we only need to choose $f$ so that the point $\<x_{n}-f\>$ of $A_1$ does not lie on the line joining any two of our points $  \<f_i \>$. 
Assume that $|K|=q$ is finite. The union of those lines has size at most ${\binom{r+1}{2}}(q-1) +r+1$.  There are  $q^{n-1}$ points 
   $\<x_{n}-f\>$ as $f$ varies.  Since we have assumed that 
   $q>(r+1)^2/2$, it follows that 
   $q^{n-1}>{\binom{r+1}{2}}(q-1) +r+1$  and  a suitable $f$ exists.
   When $K$ is infinite the argument is even easier.
    \qed

\medskip
A variant of the previous  argument can   be used 
 in characteristic 2:

  \begin {theorem} 
   \label{powerpoint independence char 2}
Let $K= \GF(2^m)$ and $d= 2^i+1$   with $(i,m)=1 $ 
and $m\ge 3$. 
Then any $4$  powerpoints  of $A_d$
are independent.
   \end {theorem}

 \proof
If $n=2$ and $d=s+1$, then
each powerpoint is spanned by $x_1^d$ or  
$(x_2+t x_1) ^d =x_2^d+tx_2^s x_1+ t^sx_2x_1^s+t^dx_1^d $
for some $t\in K$. 
By \cite[Lemma~21.3.14]{Hi}, 
the points $\<(1,0,0,0)\>$ and $\<(1, t, t^s, t^{s+1})  \>$,
$t\in K$,  form a 
$4$-independent set.  (N.\hspace{1pt}B. --
By contrast, in odd characteristic $p$, using $s=p^i$ the analogous set of points  always has 4 dependent members, so that the analogue of  the theorem does not hold.)

Now suppose that $n>2$. We are given $4$ distinct powerpoints  
 $\<f_1 ^d\>, \dots , \<f_{4} ^d\>$, and we will assume a linear dependence relation $\sum_1^{4}k_if_i^d=0$  for
  scalars $k_i$.
 Apply  the endomorphism $T$ of $A$ fixing 
 $x_1,\dots,x_{n-1}$ and sending  $x_n$ to an arbitrarily chosen  linear combination $f $ of 
 $x_1,\dots,x_{n-1}$ 
   (cf. Remark~\ref{extends to endomorphism}) 
  in order to obtain
 an identity $\sum_1^{4}k_iT(f_i)^d=0$ in the ring
  $ K [x_1,\dots,x_{n-1}]$.  If the powerpoints 
$\<T(f_i)^d\>$ are distinct then we will have reduced the number of indeterminates $x_i$, as desired: the $k_i$ are all 0.

As in the proof of   the preceding theorem, 
  we only need to {choose} $f$ so that the point  
   $\ker T = \<x_{n}-f\>$
    does not lie on the line joining any two of the
     points $\<T(f_i)\>$. 
 The union of those lines has 
 size at most ${\binom{4}{2}}(q-1) + 4$, where $q=2^m$.  There are  $q^{n-1}$ points 
   $\<x_{n}-f\>$ as $f$ varies.  
   Then a suitable $f$ exists since $m \ge 3$  implies  that 
   $q^{n-1}>6(q-1) + 4$.
   \qed
        
    \medskip    
We emphasize that the preceding theorem is a higher-dimensional generalization of a standard result in $\PG(3,q)\,$ \cite[Lemma~21.3.14]{Hi}.  In fact, since this is only a question of four points an approach that is easier than the above simply plays with
 the space spanned by the $f_i$.  
  As in Remark~\ref{code} there is an associated code  that 
  may be  worth some study. For example, an elementary examination
   of possible  dependence relations among the polynomials $f_i^d$ shows that all minimum weight codewords arise from 2-spaces of $K^n$.

\section{Independence of  power subspaces}
\label{Independence of  power subspaces}

Let $A$ be as in (\ref{eq:graded}).
Recall that, if $T$ is a subspace of $A_1$, then  $\langle T^d\rangle$ is the subspace of $A_d$ spanned by all $d$-fold products of linear polynomials in $T$. 
\vspace{2pt}
Note that
$\dim \langle T^d\rangle=  \binom{d+ \dim T -1}{d} $
since the monomials
 of degree $d$ in a basis of $T$ form a basis of $\<T^d\>$.

Before we can prove Theorem~\ref{independence} we need a few algebraic preliminaries.
In this section we will use the uncommon notation $V^{(d)}$ to denote a cartesian power, in order to distinguish it from powers in 
  rings.

%\newpage
\subsection{The universal nature of symmetric tensors} 
\label{The universal nature of symmetric tensors} 
For a   $K$-vector space $V$ and a commutative $K$-algebra $B$,
we will need an almost-basic property of symmetric $d$-multilinear $K$-forms $V{}^{(d)}\to B$; that is,  
multilinear forms $f(v_1,\dots,v_d)$ 
assuming values in $B$
and invariant under all permutations of the $v_i\in V$.

 As in Section~\ref{Background}, we view 
 the algebra $S(V)$ of symmetric tensors  
 as the polynomial algebra $A = K[X]$ for a basis $X$ of $V$,
 viewing
 $V$ as $A_1$
 and
the subspace $S(V) _d$ spanned by the $d$-fold symmetric tensors 
as  $A_d$.    A standard and elementary  universal property of symmetric tensors is the case $B=K$ of the following  

\begin{theorem}\label{B-forms}  Let $f\colon \!V^{(d)}\rightarrow B$ be a symmetric $d$-multilinear $K$-form  with values in a commutative $K$-algebra $B$ without zero divisors.  Then there is a $K$-linear mapping  $\bar{f}\colon  A_d=S(V)_d\rightarrow B$ such that$,$ for every  
 $(v_1,\ldots , v_d)\in V^{(d)} , $
$$ f(v_1, \ldots , v_d) = \bar{f}(v_1\cdot v_2\cdots v_d).$$
\end{theorem}

\proof   Let $K'$ be the field of fractions of $B$.
Then  $K'\otimes_K A=K'[X]$
and ${K'\otimes_K A_d}=  K'[X]_d$.  

Let $V'={K'\otimes_K V}$.
There is a symmetric multilinear  $K'$-form $f'$ 
determined by $f$ together with a  $K$-basis $X:= \{x_1, \ldots , x_n\}$
of  $A_1$:
for  $v'_i = \sum_j \beta_{ij}x_j \in V'$, $i = 1,\ldots ,d$, 
$\b_{ij}\in K',$  define 
\vspace{-4pt}
$$
f'(v'_1, \ldots , v'_d):= \sum_{\sigma}\Big[\prod_{i = 1}^d  \beta_{i\sigma(i) } \Big]f(x_{\sigma(1)},\ldots , x_{\sigma(d)}) ,
$$
where the above sum is over all sequences 
$\sigma=(\sigma(1),\dots, \sigma(d))$  with entries in 
$\{1,\ldots ,n\}$.   
This definition is forced by 
  multilinearity  together with  the requirement that  $f'(v'_1, \ldots , v'_d)= f (v'_1, \ldots , v'_d)$ if all $v_i'\in V$.

  By  
the field case of the theorem there is a $K'$-linear mapping $\bar{f}'\colon K'\otimes_K A_d= K'[X]_d \rightarrow K'$ such that,
for all $v_i'\in V'$, 
$$\bar{f}'(v_1'v_2'\cdots v_d') = f'(v_1', v_2', \ldots , v_d').$$
 If all $v_i'=v_i\in V=A_1$ then 
 \vspace{.5pt}
  the right side 
 %of (\ref{f'})
is just 
$ f (v_1, v_2, \ldots , v_d) \in B $.  
Hence, the desired $K$-linear mapping is $\bar f: = \bar f'|_{A_d} \colon\! A_d= K [X]_d\rightarrow B$, since 
$ f(v_1,\ldots , v_d) = f'(v_1,\ldots , v_d) = \bar{f}'(v_1\cdots v_n)$
for every 
$(v_1,\ldots ,v_d)\in V^{(d)}$.  \qed

\subsection{Proof of Theorem~\ref{independence}}

\label{Proof of Theorem independence}

We begin with the analogue of Lemma~\ref{almost independence}:

\begin{lemma}
\label{almost independence of powers}
Suppose that  $1\le r \le d$ with 
 $d!/(d-r)!\ne0 $ in $K$. 
 If $T_0$ is a subspace of $A_1$ intersecting each of $r$
  subspaces $T_1, \dots,T_r$ of $\,A_1  $ only in $0,$
then $\<T_0^d\>\cap\<T_1^d, \dots ,T_r^d \>=0.$
\end{lemma}

\proof 
  It suffices to prove that  there is a subspace $N_{0}$ of $A_d$  such that $N_{0}$ contains $T_j^d$ for all $j\ge 1$ and $N_{0}\cap \<T_{0}^d\> = 0$.    
Change coordinates in $A_1$ so that  $x_1, ... ,x_k$  is a basis of 
$T_{0}$.  Let $B:=K[x_1, ... ,x_k]$, so that  
 $B_d:=B\cap A_d$ is $\<T_0^d\>$.

If $1\le j\le r$ then  $T_0\oplus T_{j}$ is a direct summand of $A_1$, so that there is a linear transformation $\lambda_j \colon\! A_1 \to A_1$ such that $\lambda_j (x_i) = x_i, \ i= 1,\ldots , k$, and $\lambda_j(T_j) = 0$.   (The behavior of $\lambda_j$ on a complement to 
$T_0\oplus T_{j}$ in $A_1$ is irrelevant to the proof.)

  If $r<j\le d$ then 
$\lambda_j \colon\! A_1 \to A_1$ will be the identity map.

Let $\Theta$ be a (left) transversal for the pointwise stabilizer 
$S_{d-r}$ of 
$1,\dots ,r$ in the symmetric group $S_d$ on $\{1,\dots,d\}$, so that $| \Theta|=d!/(d-r)!$.  
 Define a $d$-multilinear $K$-form 
$ L_{0} \colon\! A_1^{(d)} \to A_d$
by \vspace{-2pt}
$$
%\vspace{-6pt}
%\begin{equation}\label{eq:form}
 L_{0}(v_1, \ldots  ,v_{d}):= \sum_{\pi\in \Theta} \prod_{j=1}^d\lambda_j(v_{\pi(j)}) .
$$
%\vspace{-2pt}
%\end{equation}
We claim that {\em $L_0$ is symmetric.}  For, let 
$\pi\in\Theta,$
$\rho\in S_d$, and write $\rho\pi=\pi'\sigma$ with $\pi'\in\Theta,$ $ \sigma\in S_{d-r}$.
Then
$$
\prod_{j=1}^r \lambda_j(v_{\rho\pi(j)}) \!\!
\prod_{j=r+1}^d \lambda_j( v_{\rho\pi(j)})
\!=\!
\prod_{j=1}^r \lambda_j(v_{\pi'\sigma(j)})\!\!
\prod_{j=r+1}^d   v_{\pi'\sigma(j)}
\!=\!
\prod_{j=1}^r \lambda_j(v_{\pi'(j)})\!\!
\prod_{j=r+1}^d   v_{\pi'(j)} ,
$$
since $\{{\pi'\sigma (j)}\mid r+1\le j\le d\}$ is the complement
in %the full set
 $\{1,\dots,d\}$ %of indices   
 of 
$\{{\pi'\sigma(j)}\mid 1\le j\le r\}
\break
=\{{\pi'(j)}\mid 1\le j\le r\}$, and hence is $\{{\pi'(j)}\mid r+1\le j\le d\}$.
Consequently, $\rho$ permutes the summands that define $L_0$, which proves the claim.

By   Theorem~\ref{B-forms},    there is a linear transformation
$\bar L_{0}  \colon\! A_d \to B$
such that
\vspace{-1pt}
\begin{equation}
\label{eq:form}
\bar L_{0}(v_1 \cdots v_d ) =L_{0}(v_1,\ldots , v_d)
 = \sum_{\pi\in \Theta} \prod_{j=1}^d\lambda_j(v_{\pi(j)}) 
\end{equation}
for all ${v_i\in A_1}$.
  Clearly $\bar L_{0}(A_d)\subseteq A_d$.
We will show that {\em  $N_{0} := \ker\bar L_{0}$  behaves as required  at the start of this proof}.

Consider $j\ge 1$.  If  all $v_i\in T_j$, then 
$v_{\pi(j)}\in T_j
\subseteq  \ker \lambda_j,$ and   each summand on the right side of  (\ref{eq:form}) is $0$.  Thus, 
 $\bar L_{0}(T_j^d)=0$.

It remains to determine 
%\marginpar{NEWWWW}
the action of 
$\bar L_{0} $ on $\<T_{0}^d\>$. We first calculate $\bar L_{0} $ on each monomial  $x^{\alpha} =  x_1^{a_1}\cdots x_k^{a_k}$,   $\sum _ia_i = d$.  Since $\lambda_j(x_i) = x_i$ for all $i\le k$ and all $j$,  (\ref{eq:form}) gives  
  \begin{equation}\label{eq:form2} 
  \bar L_{0} (x^{\alpha}) =  |\Theta| \: x_1^{a_1}  \cdots x_k^{a_k}.
   \end{equation}
   Here $|\Theta|=d!/(d-r)! \ne0 $ 
   by hypothesis.
      
   Thus, as $x^{\alpha}$ ranges over all monomials of degree $d$ in $x_1, ... ,x_k$, their $\bar L_{0} $-images form a $K$-basis for $B_d$.  Consequently, $\bar L_{0} $  restricted to $\<T_{0}^d\>$  
 is a   surjection $ \langle T_{0}^d\rangle\rightarrow B_d$,
 and hence an 
 isomorphism since   $\dim\langle T_{0}^d\rangle=\dim B_d$.  Thus, 
 $N_{0}\cap \langle T_{0}^d\>=\ker \bar L_{0} \cap \langle T_{0}^d\rangle = 0$,  as required.  \qed
 \medskip
  
\paragraph {\bf Proof of Theorem~\ref{independence}}
 The case $e=1$ of  Theorem~\ref{independence} follows immediately
  from the preceding lemma.  The general case is obtained exactly as in  the proof of {\rm Theorem~\ref{Veronesean e+1 independence}}
 near the end of  Section~\ref{Families of subspaces}.  \qed

\medskip
 When $r=d$, an entirely different proof of Theorem~\ref{independence} is  
 obtained by combining Theorems~\ref{Veronesean independence} and \ref{th:pover}.    
    Theorem~\ref{independence} 
clearly contains Theorem~\ref{power point independence} 
 as a special case, but it     does not quite contain
   Theorem~\ref{powerpoint independence  small char}:  the requirements on $r$ are less stringent in the latter result.  (For example, if $d=5$  and the characteristic is $r=3$, then $3$ divides 
   $5!/(5-3)!$ but none of the binomial coefficients~$\binom{5}{i}$.)%

     \section{$r-$Independence of  power subspaces}
\label{3-Independence of  power subspaces}
In this section we will use subspaces of polynomials to prove a (weak) variation on the results in the preceding section:

\begin{theorem}
\label{3-independence}

Let $r\ge1$. 
If  $d>1$ is not a power of  $ \Char K $ and if 
  $\mathcal T$ is any $r-$independent set of subspaces of $A_1,$   then  $\{ \<T ^d\>\mid T\in {\mathcal T  }\}$ is an $ r+1 -$independent
 set  in~$A_d$.

 \end{theorem}

	\subsection{Calculating with spaces of polynomials}
\label{Calculating with subspaces of polynomials} 
 We will make frequent use of the following elementary observation and its consequences.
 
\begin{proposition}\label{th:dim} 
In {\rm(\ref{eq:graded})} let $U_1$ and $U_2$ be subspaces of $A_1$ such that $U_1\cap U_2 = 0$.   If $d>1,$ then
\begin{eqnarray}
    \langle A_{d-1}U_1\rangle\cap\langle A_{d-1}U_2\rangle
    \hspace{-6pt} &=&\hspace{-6pt} 
    \langle A_{d-2}U_1U_2\rangle
    \label{eq:seven}
    \\
 \langle (U_1 + U_2) ^d\rangle \hspace{-6pt}&=&\hspace{-6pt} \bigoplus_{k=0}^{d} \langle U_1^kU_2^{d-k}\rangle 
 \label{eq:eight} \\
\< U_1^d\>\cap \langle A_{d-1}U_2\rangle \hspace{-6pt}&=& \hspace{-6pt}0.
\label{eq:nine}
\end{eqnarray}

\end{proposition} 

\proof 
 Let $X_1 := \{x_1, \ldots , x_{\ell}\}$ and $X_2:= \{ x_{\ell + 1},\ldots , x_m\}$ be respective  bases for   $U_1$ and $U_2$.  
 Let  $ X\supseteq X_1\dot\cup X_2$ be a basis 
  of $V$.
 
  Both sides of each  of the above equations are subspaces of $A_d$.  The left side of (\ref{eq:eight}) is the subspace of $A_d$ spanned by all monomials of degree $d$ with factors chosen from $X_1 \dot\cup X_2$.  Partitioning these monomials according to the number of factors of $X_1$ they contain proves  (\ref{eq:eight}).

For  (\ref{eq:nine}) note that $ \<U_1^d\>$
is spanned by  monomials in $X_1$ of degree $d$, while $\langle A_{d-1}U_2\rangle$ is spanned by   monomials   containing at least one factor from~$X_2$. 

For (\ref{eq:seven}),   
 consider the following pairwise disjoint  sets
of  monomials in $X$:
\begin{itemize}
%  \item $Y_0:=M_d(X_0)$ is the set of monomials
 %  of $M_d(X)$ with no factor in $X_1 \dot\cup X_2$,
  \item  $Y_i\subset A_{d-1}X_i $  is the set of monomials in $X$ with at least one factor from $X_i$  and no factor  from $X_{3-i}$, for $i=1,2$, and
  \item   $Y_{12}\subset A_{d-2}X_1X_1$ is the set of all monomials in $X$ having at least one factor from $X_1$ and at least one from $X_2$.  
\end{itemize}  
  It follows that  
  $\langle Y_1\rangle\cap \langle Y_2\rangle =0$    and   
$\langle A_{d-1}U_i\rangle = \langle Y_i\rangle\oplus\langle Y_{12}\rangle  $ for $ i = 1,2$.
 Consequently, $\langle A_{d-1}U_1\rangle \cap \langle A_{d-1}U_2\rangle = \langle Y_{12}\rangle  = 
    \langle A_{d-2}U_1U_2\rangle $.~\qed
  \medskip

We can now show   that the $d$th power operator   commutes with intersections: 
\begin{corollary}\label{cor:ten}  For any subspaces $B$ and $C$   of $A_1 ,$  
\begin{equation}\label{eq:d-cap}
\<B^d \>\cap\< C^d \>= \<(B \cap C)^d \>.
\end{equation}
\end{corollary}

\proof We may  assume that $d>1$.  Set
 $C_1: = B\cap C$, and choose a subspace $C_2$ such that $C = C_1\oplus C_2$.  Since $d>1$,   (\ref{eq:eight}) yields  
 \begin{equation}\label{eq:eleven}
  \<C^d \rangle = \bigoplus_{j=0}^{d} \langle C_1^jC_2^{d-j}\rangle=  \<C_1^d \rangle\oplus \langle C_2C^{d-1}\rangle.
\end{equation}    

Since $B\cap C_2=  B\cap (C\cap C_2)= C_1\cap C_2  = 0$,   
(\ref{eq:nine}) forces 
$\langle B^d \rangle\cap \langle C^{d-1}C_2\rangle = 0$.  On the other hand, $\langle B^d \rangle $ contains $\langle C_1^d \rangle $, the first summand at the end of (\ref{eq:eleven}).  Thus, 
$$
\<B^d \rangle\cap \langle C^d \rangle  
= \<C_1^d \rangle
= \<(B\cap C)^d \rangle. \ \ \ \ \qed
$$   
%\smallskip 
 
  \subsection{Proof of Theorem~\ref{3-independence}}
 \label{Proof of Theorem 3-independence}

We begin with a special case:

 \begin{proposition}
\label{k to k+1}
Suppose that $A_1 = T_1\oplus \cdots \oplus T_r$ with $\dim T_i = s,$
and that 
$T_{r+1}$ is an $s$-space in $A_1$ 
such that  the set $\{T_1,\dots,T_{r+1} \}$ is $r$-independent.   Then 
$\langle T_1^d, \dots, T_r^d \rangle=
\langle T_1^d \rangle\oplus \cdots\oplus\langle T_r^d \rangle,$
and one of the following holds$:$
\begin{enumerate} 
  \item $\langle T_{r+1}^d \rangle\cap (\langle T_1^d \rangle\oplus \cdots \oplus\langle T_r^d \rangle) =   0 ,$ or 
  \item  ${\rm char} K $ is a  prime $p,$    $d$ is a power of $p,$  and  
 $ \dim [\langle T_{r+1}^d \rangle\cap 
 (\langle T_1^d \rangle\oplus \cdots \oplus\langle T_r^d \rangle)]
    = s.
 $
  \end{enumerate}
\end{proposition}
%\newpage

\proof  Let  $\{y_1, \ldots , y_s\}$  be a basis of $T_{r+1}$.
If $x_{ij}$ is the projection of $y_j$ into $T_i$, $1\le i\le r$, then
$y_j=\sum_ix_{ij}$ with
each $x_{ij}\ne 0$ due to $r$-independence, and  
$X_{i}:=\{x_{ij}\mid 1\le j\le s\} $ is a basis of $T_i$ for $i\le r$.
%  Let  $X= \dot \cup _iX_i$. 
Then 
$\langle T_1^d ,\dots ,  T_r^d \rangle=
\langle T_1^d \rangle\oplus \cdots \oplus\langle T_r^d \rangle
$
since 
$\langle T_i^d \rangle $ is  spanned   by monomials in $X_i$.

If $0\ne f\in \langle T_{r+1}^d \rangle\cap 
(\langle T_1^d \rangle\oplus \cdots \oplus\langle T_r^d \rangle)$,      then
\begin{equation}\label{eq:t3d} 
 f = \sum_{\a}
k_\a  \prod_{j=1}^s y_j^{a_j},
 \end{equation}
 %\overset{\scriptstyle m=(a_1,\ldots , a_k)}
 %{\in M[1,k]_d}  WMK
 where $ k_\a\in K$ and the sum is indexed by all
 $\a=(a_1,\ldots , a_s) $ with all $a_i\ge0$ and $\sum_ ia_i=d$.
  Since $f\in   \langle T_1^d \rangle\oplus \cdots \oplus\langle T_r^d \rangle $,  when expanded as a linear combination of monomials in 
  $\cup_iX_i$ of degree $d$ the coefficients of monomials in (\ref{eq:t3d}) with ``mixed terms'' -- i.\:e., monomials containing members of  $X_iX_j$
  with $i\ne j$ -- must be zero.
 
 Let $\a = (a_1,\ldots ,a_s)$ be as above and suppose that (at least) two of the numbers $a_ \ell $ and $a_{m}$ are positive ($\ell\neq m$).  Then 
 the product  $\displaystyle \prod_{j=1}^s x_{{* j}}^{a_{j}}$, 
where  $*=1$   except that $*=2$ when $j=m$,
contains a term in $X_\ell X_m$, 
and this product occurs  only once in (\ref{eq:t3d}).
 Since $f\in   \langle T_1^d \rangle\oplus \cdots \oplus\langle T_r^d \rangle $, 
it follows that the coefficient $k_\alpha$ in  (\ref{eq:t3d}) is zero for all  $\a = (a_1,\ldots , a_s)$ having at least two      nonzero  terms.

Now   (\ref{eq:t3d}) reduces to  
\begin{equation}\label{eq:t3d2}
f = \sum_i k_i y_i^d.
\end{equation}
By  the Binomial Theorem, $f$ involves  a nonzero mixed term  containing a member  of  $X_1X_2$ unless $K$ has characteristic $p>0$ and $d=p^e$ for some $e$. Then $ y_i^d = (\sum_ix_{ij})^d = \sum_ix_{ij}^d 
\in \langle T_1^d \rangle\oplus \cdots \oplus\langle T_r^d \rangle,$ and   (2) holds.~~\qed

\medskip \medskip

%\newpage
%\vspace{.15in}

\noindent  \emph{Proof of} Theorem~\ref{3-independence}.   
 \emph{It suffices to prove that$,$  if $\,T_1,\ldots ,  T_{r+1}$ are distinct members 
of $\mathcal H,$ then}
\begin{equation} \label{k-ind-var}
\langle T_{r+1}^d\rangle \cap\langle T_1^d,\ldots , T_r^d\rangle = 0 
\end{equation} 
(compare Lemma~\ref{almost independence of powers}).
Set $ T:= \langle T_1,\ldots , T_r\rangle = 
\langle T_1 \rangle\oplus \cdots \oplus\langle T_r \rangle\,$
(by $r$-independence) 
with corresponding projections
 $\pi_i\colon T\rightarrow T_i$.

Let
 $U_{r+1} := T_{r+1} \cap T$.   
We may assume that $U_{r+1}\neq 0$, as otherwise
  $\langle T_1,\ldots ,   T_{r+1} \rangle\break
=T_1\oplus \cdots \oplus T_{r+1} $ and  
(\ref{k-ind-var}) is clear (use  Corollary~\ref{cor:ten}
with $C=T$).

Once again let $\{y_1, \ldots , y_s\}$ be a basis for $U_{r+1}$  
and $x_{ij}: = \pi_i(y_j)$, $1\leq j\leq s$, $1\leq i\leq r$.  If, for some $i$, $\{ x_{ij} \mid 1\leq j\leq s\}$ is   linearly dependent, then some  $y\ne0$ in $U_{r+1}$  satisfies $\pi_i(y) = 0$,  and
 $r$-independence produces  the contradiction 
 $y\in U_{r+1}\cap \ker\pi_i\leq T_{r+1}\cap \oplus_{j\ne i}T_j=0$.

 Thus,   $U_i:=\< x_{ij} \mid 1\leq j\leq s\>$ is an $s$-subspace  of $T_i$  for each $i$.
Since $y_j\in T$ we have 
$y_j\in \< \pi_i(y_j) \mid 1\leq i\leq r\>$ 
and hence $U_{r+1}\le U := 
 \<U_1,\dots, U_{r}\>$.
Since   $U_i\le T_i$ and $\{T_1,\dots, T_{r+1}\}$
 is an $r$-independent family, so is the family 
  $\{U_1,\dots, U_{r+1}\}$ of subspaces of $T$.  
Apply Proposition~\ref{k to k+1} to this family
 with $  U_i$ and $U $ in the roles of $  T_i
$ and $ A_1$: 
\begin{equation} 
\label{U's}
\< U_{r+1}^d\rangle\cap \<U_1^d,\dots, U_{r}^d\> = 0,
\end{equation}
 which resembles our goal (\ref{k-ind-var}).

We claim that \begin{equation}\label{eq:uk+1}
\langle T_{r+1}^d\rangle \cap\langle T_1^d,\ldots , T_r^d\rangle 
\le\langle U_{r+1}^d\rangle \cap\langle T_1^d,\ldots , T_r^d\rangle  
\end{equation}
(and later we will show that the right hand side is 0). 
For, select a complement $W_i$ to $U_i$ in $T_i$ for $i = 1,\ldots , r+1$,
and let $W:= \langle W_1,\ldots , W_r\rangle$.
Then  $T = U\oplus W$
and $
\langle T_1,\ldots , T_{r+1}\rangle = U\oplus  W \oplus W_{r+1}.
$
By Corollary~\ref{cor:ten},  $\langle T_{r+1}^d\rangle \cap 
\langle T_1^d,\ldots , T_r^d\rangle
\leq\langle T_{r+1}^d\rangle \cap \langle T^d \rangle= \langle (T_{r+1}\cap T)^d\rangle = \langle U_{r+1}^d\rangle$,
which proves (\ref{eq:uk+1}).

By  (\ref{eq:nine}), $\<U^d\rangle \cap \langle A_{d-1}W\rangle = 0 $.
Since $\<U_1^d,\dots, U_{r}^d\> \leq \langle U^d\rangle$, 
by the modular law   
$$
\langle U^d\rangle \cap [ \<U_1^d,\dots, U_{r}^d\>) \oplus \langle A_{d-1}W\rangle ] \! = \! \<U_1^d,\dots, U_{r}^d\> \oplus (\langle U^d\rangle \cap \langle A_{d-1}W\rangle ) \! = \! 
\<U_1^d,\dots, U_{r}^d\> ,
$$
and hence  (since $T_i^d=(U_i\oplus W_i)^d\le U_i^d\oplus A_{d-1}W_i   $
by (\ref{eq:eight}))
\begin{equation*}
\label{int}
\begin{array} {lllll}
\langle U^d\rangle \cap \langle T_1^d,\ldots , T_r^d\rangle   
\hspace{-6pt}&\le&\hspace{-6pt}\langle U^d\rangle \cap [\<U_1^d,\dots, U_{r}^d\>  + \langle A_{d-1}W_1,\ldots , A_{d-1}W_k\rangle]  \\
\hspace{-6pt}&=&\hspace{-6pt}\langle U^d\rangle\cap [\<U_1^d,\dots, U_{r}^d\>  + \langle A_{d-1}W\rangle ] \\
\hspace{-6pt}&=&\hspace{-6pt} 
\<U_1^d,\dots, U_{r}^d\> .
\end{array}
\end{equation*}
Since $\langle U_{r+1}^d\rangle\leq \langle U^d\rangle$,  
%(\ref{int})
(\ref{U's}) yields 
\begin{eqnarray*}
\langle U_{r+1}^d\rangle\cap 
\langle T_1^d,\ldots , T_r^d\rangle  
\hspace{-6pt}&=&\hspace{-6pt} 
\langle U_{r+1}^d\rangle\cap [\langle U^d\rangle \cap 
\langle T_1^d,\ldots , T_r^d\rangle 
]\label{uk11} \\
\hspace{-6pt}&\le &\hspace{-6pt} 
\langle U_{r+1}^d\rangle\cap 
\<U_1^d,\dots, U_{r}^d\>
 =0 .
 \end{eqnarray*}
Now    (\ref{eq:uk+1}) implies
 (\ref{k-ind-var}), as required. 
  \qed  
  
 %\newpage

 \section{Generalized dual arcs and polynomial rings} 
 \label{Generalized dual arcs and polynomial rings}    
\subsection{\hspace{-6pt} Definitions} \hspace{-6pt}
A {\em generalized dual arc of   $V=K^n$
 with  vector dimensions 
 $(n, n_1,  \ldots , n_d)$} is a set ${\mathcal D}$ of $n_1$-subspaces of   $V$ such that the intersection of any $j$ of them has dimension $n_j>0$, $j = 2,\ldots ,d $, and the intersection of any $d+1$
  of them is 0.  
 This notion was introduced in \cite{KSS1, KSS2}, but using projective dimension instead of vector space dimension. 
  When  $d = 2$ and $n_d = 1$, ${\mathcal D}$ is   a {\em  dual arc}.

 Suppose ${\mathcal D}$ is a generalized dual arc   with  vector dimensions   $(n, n_1,  \ldots , n_d)$.   If $D\in {\mathcal D}$, then 
  ${\mathcal D}':= \{ D\cap D' \mid D'\in {\mathcal D}-\{ D\}\}$ is a generalized dual arc with   vector dimensions    
  $( n_1,  \ldots , n_d)$. In general, this procedure can be iterated.

 For any subset $\Delta$ of a sublattice ${\mathbf L}$ of
 the lattice
 ${\mathbf L}(V)$  of all subspaces of $V$, let ${\mathbf L}(\Delta)$ denote the sublattice generated by $\Delta$ (the smallest sublattice containing $\Delta$),
 and  
  let ${\mathbf I}({\Delta})$ denote the {\em ideal} generated by $\Delta$ (the set of all elements of  ${\mathbf L}(\Delta)$ that   are bounded above by at least one element of $\Delta$).  
We call   $\Delta$ {\em regular} if  $V=\langle \Delta \rangle $, and, for each   intersection $U$ of finitely many members of $\Delta $, 
 \begin{equation}\label{eq:reg-lattice}
 U = \langle  U\cap D \mid D\in \Delta,~ D\not\subseteq U  \rangle.
 \end{equation}
We call $\Delta$  {\em strongly regular} if 
 it is regular and 
\begin{equation}\label{eq:sr-lattice}
  U\cap  \langle D_1,   \dots  , D_\ell \rangle =
\langle   U\cap   D_i  \mid 1\le i\le \ell\rangle
\end{equation}
for any $D_1,   \dots  , D_\ell\in \Delta  $ and any subspace $U$ that is an intersection  of finitely many members of  $\Delta$.  There are many   stronger versions of this concept possible\footnote{For example,  $U$ might be restricted to 
range over all  elements of  ${\mathbf L}(\Delta)\cap{\mathbf I}(\Delta)$.}, but this definition is geared to conform to the definitions appearing in \cite{KSS1, KSS2}. 

\subsection{An elementary  construction}
\label{A construction}
Consider (\ref{eq:graded}), fix $d\geq 2$, and let ${\mathcal D}$ be the following set of   $K$-subspaces of $A_d$:
   \begin{equation}\label{generalized dual arc}
    {\mathcal D}:=\{ A_{d-1} y  \mid  0\ne y\in A_1  \}.
     \end{equation}
 If $z = \alpha y,   \alpha\in K^*$, then $A_{d-1}y = A_{d-1}z$, so ${\mathcal D}$ is parameterized by the $1$-spaces of $A_1$.  
If $D_j: = A_{d-1}y_j$, $1\le j\leq d$, are distinct elements of ${\mathcal D}$,  {\em we claim that}
    \begin{equation}\label{eq:gda}
     D_1\cap D_2\cap \cdots \cap D_j = A_{d-j}y_1\cdots y_j.  
     \end{equation}
    Clearly the right side of   (\ref{eq:gda}) is contained in the left. The left  side is precisely  all homogeneous polynomials of degree $d$ that are divisible by all of the  polynomials $y_1,..., y_j$.  Since the polynomials $y_i$ are      primes
    of the  unique factorization  domain $A$
no two 
of which are associates, each   polynomial  on the left side of (\ref{eq:gda}) is a multiple of  $y_1\cdots y_j$ and so lies in the right side, as claimed.
    
    The  dimension of the subspace   in   (\ref{eq:gda}) is $\dim A_{d-j}$.  This proves 
    
    \begin{theorem}
    \label{dual arc construction}
    ${\mathcal D}$   is a generalized dual arc in $A_d$  with vector dimensions
    $$\mbox{$\Big(\binom{n+d-1}{d}, \binom{n+d-2}{d-1}, \ldots , n_{d-1} = \binom{n}{1}=n, n_{d} = 1 \Big).$} $$
    \end{theorem}
    
  In  \cite{KSS2} there are objects with similar properties constructed without the use of polynomials.
  Note that  $\< {\mathcal D} \>=A_d$.
%\newpage
\subsection{Regularity properties of this construction: Examples.}\

\Example \rm {\em Strong regularity can fail for} (\ref{generalized dual arc}).
 Suppose $n= \dim A_1 \geq 3$ and $d=2$. 
 We will see in the next Example that  ${\mathcal D}$
 in (\ref{generalized dual arc}) is regular since $d$ is small.
 However, ${\mathcal D}$ {\em is not strongly regular}.  
For, consider the subspaces $D_i: = A_1x_i,$ $ i = 1,2$, and
 $U: = A_1(x_1+x_2)< \langle D_1, D_2\rangle$ belonging to  
 $ {\mathcal D}$.

 Since $U \cap  D_i = \<x_i(x_1+x_2)\>$, $ i=1,2$, 
$  \dim \langle U\cap D_1 , U\cap D_2 \rangle=2<n
= \dim U =
\dim (U\cap\langle D_1, D_2\rangle ) ,$
so that $ U\cap\langle D_1, D_2\rangle   \neq \langle U\cap D_1 , D_2\cap D'\rangle  $ 
and   (\ref{eq:sr-lattice})  fails, as required.

%\vspc
\Example \rm
{\em Even regularity can fail for} (\ref{generalized dual arc}),  {\em but only for finite  $K$ and enormous   $d$.}
Since $\< {\mathcal D} \>=A_d$,  in order to see this we must examine 
(\ref{eq:reg-lattice}).  

Each $U\ne0$  in (\ref{eq:reg-lattice}) has the form
$U: =\bigcap_{i=1}^k ( A_{d-1} y_i)
=A_{d-k}\pi$ for distinct $ \<y_1 \>,\dots , \<y_k\>$
in $A_1$ and  $\pi:=y_1\cdots y_k$,  where  $d\ge k$.
Regularity asserts  that each such $U$  
 is spanned by the subspaces $U\cap (A_{d-1}y)$ where $y$ ranges over the set $Y'$ of linear polynomials 
 {\em not} in
    $  \<y_1 \> \cup \dots \cup \<y_k\> $.  
		
By (\ref{eq:gda}),  
\begin{eqnarray*} 
 \langle U\cap (A_{d-1}y) \mid y\in Y' \rangle \hspace{-6pt} &=&\hspace{-6pt}  \langle A_{d-k-1}y\pi \mid y\in Y'\rangle \\
  &=&\hspace{-6pt}  \langle A_{d-k-1}Y' \rangle\pi .  \end{eqnarray*}
 In particular, if   $A_1 = \langle  Y'\rangle$  then 
 $ \langle U\cap (A_{d-1}y) \mid y\in Y' \rangle=A_{d-k}\pi = U$.

Thus, if ${\mathcal D}$ is {\em not} regular  then,
for some $U$, the corresponding set  $Y'$ must span a proper subspace $H$ of $A_1$.  Then  the $k$-set  $ \{\<y_1\>,\ldots, \< y_k\>\}$ must contain  all $1$-spaces of $A_1$ not in $H$. 
Thus,  $K  $ is  finite  and $k\geq |K|^{n-1}$, the number of points not in a hyperplane of $\P(A_1)$.  Then $d \geq k \geq |K|^{n-1}$.

Conversely, if $d \ge  |K|^{n-1}$ choose
$U: =\bigcap_{i=1}^k ( A_{d-1} y_i),$ where
$ \{\<y_1\>,\ldots, \< y_k\>\}$   consists of all $k=|K|^{n-1}$  points outside the hyperplane $H=\<x_1,\dots,x_{n-1}\>$  of $A_1$.  Then the previous argument produces a subspace 
$$ \langle U\cap (A_{d-1} y) \mid y\in Y' \rangle = \langle A_{d-k-1}H\rangle \pi
$$
of dimension smaller than that of $U = A_{d-k}\pi$, which proves nonregularity.   

 Thus, we have proved 
\begin{proposition}\label{th:non-reg}  A generalized dual arc  ${\mathcal D}$  in {\rm (\ref{generalized dual arc})}  need not be strongly regular.  

Moreover$,$  ${\mathcal D}$ is not regular 
if and only if 
$K$ is a finite field $ \GF(q)$  and $d \geq q^{n-1}$.
\end{proposition} 

\subsection{A generalization of the construction}  Let $I_k$ 
denote  the set of $1$-spaces spanned by   the various homogeneous polynomials of degree $k$ that are powers of irreducible polynomials
(the degrees of these irreducible polynomials are allowed to vary).  If $d\ge k$,  then 
$$ {\mathcal D}: = \{ A_{d-k}y  \mid y\in I_k\}
$$
 is a set of $ |I_k|$  subspaces of $A_d$, each  of dimension $\dim A_{d-k}$.  Choose the positive integer $c$ such that $0\le d-kc<k$.  By unique factorization in 
 $A$, if $2\leq m\leq c$  then the intersection of any $m$ members of ${\mathcal D}$ is a subspace of dimension $\dim A_{d-mk}$, and the intersection of any $c+1$ members of ${\mathcal D}$ is $0$.  
 Thus,
 
 \begin{theorem}
 ${\mathcal D}$  is a generalized dual arc of   vector dimensions
$$ (\dim A_d, \dim A_{d-k}, \ldots , \dim A_{d-ck}).$$
 \end{theorem}

Theorem~\ref{dual arc construction}   is the special case $k=1$.
 This time,  ${\mathcal D}$ need not span $A_d$ (e.\:g., if 
 $k=\Char K$).

\subsection{Further variations}
\subsubsection{An example leading to a $3$-independent  family.} 
\begin{lemma}\label{lem:dim-ah}  If $\dim_K A_1 = n$ and  $H$ is a $k$-space in $A_1,$ then $A_1H$ is a subspace of $A_2$ of dimension $kn -\binom{k}{2}$.
\end{lemma}

\proof  Clearly 
\begin{equation} \label{eq:a1h}
\langle A_1H\rangle = A_1x_1 + A_1x_2 +\cdots + A_1x_k,
\end{equation}
assuming that $\{ x_1,\ldots , x_k\}$ is a basis of $H$.  
The intersection of any two summands on the right side of  (\ref{eq:a1h}) is a 1-space while the intersection of any three is $0$.  
Consequently, by elementary linear algebra and inclusion-exclusion, 
$$  
 \dim \big( \sum_{i=1}^k A_1{x_i} \big)= \  \sum_{i=1}^k \dim A_1x_i  - 
 \hspace{-4pt} \sum_{1\leq i<j \leq k} 
 \dim (A_1x_i\cap A_1x_j  ) %\\
=  kn - \binom{k}{2} .  \  \    \qed
  $$

  \begin{lemma}
  \label{partial spread lemma}
  Let ${\mathcal S}$ be a partial spread of $k$-subspaces of $A_1,$ where $\dim A_1 = 2k$.   
  \begin{enumerate} 
    \item If $H\in {\mathcal S}$ then $\dim  \langle A_1H\rangle  = k(3k+1)/2$.
    \item For any distinct $H_1, H_2 \in {\mathcal S}$, $\langle A_1H_1\rangle \cap \langle A_1H_2\rangle = \langle H_1H_2\rangle$ has dimension~$k^2$.
    \item For any distinct $H_1, H_2, H_3\in {\mathcal S}$, \  $\dim \big(\bigcap_{i=1}^3 \langle A_1H_i\rangle \big )=\binom{k}{2}$.
  \end{enumerate}
\end{lemma}
\proof (1)  follows from Lemma \ref{lem:dim-ah} with $n = 2k$. 

(2)  follows from (\ref{eq:seven}).  

(3) 
We may assume that $X_1=\{x_1,\dots,x_k\}$ and 
$X_2=\{x_{k+1},\dots,x_{2k}\}$ are bases of $H_1 $  and $H_2$, respectively, such that $ \{ x_i+x_{k+i}\mid 1\le i\le k\}$  
is a basis of $H_3$.   
Using (2), if $y\in \langle A_1H_1\rangle \cap \langle A_1H_2\rangle  \cap \langle A_1H_3\rangle
= \langle H_1H_2\rangle \cap \langle A_1H_3\rangle$ then we can write 
\begin{equation} \label{eq:int-123}
 y = \sum_{i=1}^k  y_i (x_i + x_{k+i})
\end{equation}
where $y_i = \sum_{j=1}^k (\alpha_{ij}x_j + \beta_{ij}x_{k+j})$.  Since  $X_1X_2$ is a basis of $\langle H_1H_2\rangle$,  
for $1\le i,j\le k$ the coefficients of $x_i^2, x_ix_j$ and $x_{k+i}x_{k+j}$   in the polynomial $y$ must be $0$:\begin{equation} \label{eq:ab-skew}
\alpha_{ii} = \beta_{ii} = 0, \alpha_{ij} = -\alpha_{ji}, \mbox{ and } \beta_{ij} = -\beta_{ji}, 1\le i,j\le k . 
\end{equation} 

The coefficient of $x_ix_{k+j}$ in  (\ref{eq:int-123}) is $\beta_{ij} + \alpha_{ji} = \beta_{ij} - \alpha_{ij}$, by (\ref{eq:ab-skew}).   Similarly, the coefficient of $x_jx_{k+i}$ 
is $ \alpha_{ij}-\beta_{ij} $.  Then
\begin{equation}\label{eq:cap123}
 y = \sum_{1\le i<j\le k} (\beta_{ij} - \alpha_{ij})(x_ix_{k+j} - x_jx_{k+i} ).
 \end{equation}

Conversely, by reversing the steps, it is clear that any $y$ 
as  in (\ref{eq:cap123}) lies in $\langle H_1H_2\rangle = \langle A_1H_1\rangle\cap \langle A_1H_2\rangle$ and also in $\langle A_1H_3\rangle$.  
It follows that the desired dimension is the number of pairs $i,j$ such that $1\le i<j\le k$.\qed

\medskip
 
   Again
consider a partial spread  ${\mathcal S}$     and ${\mathcal D } = \{ A_1H \mid H\in {\mathcal S}\}$ 
in  Lemma~\ref{partial spread lemma}:
  a set of subspaces 
 of the $\binom{2k+1}{2}$-space~$A_2$
 \vspace{2pt}
  of dimension $k(3k+1)/2$, any two meeting in a space of dimension $k^2$  and any three meeting in  a subspace of dimension  $k(k-1)/2$.  Subtracting these numbers from $\dim A_2$, we see that these intersections  have codimensions $\binom{k+1}{2}$, $2 \binom{k+1}{2}$, and $3 \binom{k+1}{2}$, respectively.   This implies the following:

\begin{proposition} 
\label{dual spread}
Let ${\mathcal S}$ be a  partial spread of  $k$-spaces of the $2k$-space $A_1,$  and consider the set 
 ${\mathcal D} = \{ A_1H \mid H\in {\mathcal S}\}$ of  subspaces of $A_2$.  The  dual set ${\mathcal D}^*$  is a $3$-independent family of $\binom{k+1}{2}$-spaces in the dual $\binom{2k+1}{2}$-space~$A_2^*$.
  \end{proposition}

Intersection dimensions of four members of the preceding set ${\mathcal D}$   are not, in general, constant.  

\subsubsection{Three $4$-space structures in dimension $10$}
Let $K=\GF(q)$.

\Example \rm  Let $n = \dim A_1 = 3$  and ${\mathcal D}_1= \{ A_2f \mid  0 \neq  f \in A_1\}$.  Then ${\mathcal D}_1$ is a set of $1+q+q^2$ \ $6$-subspaces of the $10$-space $A_3$.  Any two members of ${\mathcal D}_1$ meet at a $3$-space and so generate a $9$-space.  Since the intersection of any three is a $1$-space, any three span a space of dimension $6+6+6-3-3-3+1 = 10$, hence the entire space.  Then in the dual space $A_3^*$   
we obtain a set ${\mathcal D}_1^*$ of $1+q+q^2$ \ $4$-spaces, any two of which meet at a $1$-space, and any three of which meet at $0$.  Thus, ${\mathcal D}_1^*$ {\em is a dual arc 
in $A_3(K^{(3)})^*$ with vector dimensions $(10, 4,1)$.}

\Example \rm If $n = 4$ then ${\mathcal D}_2 :=\{ A_1f \mid 0 \neq  f \in A_1\}$  consists of $1+q+q^2 + q^3$ \ $4$-spaces of the $10$-space $A_2$.  Any two members of ${\mathcal D}_2$ intersect at a $1$-space, so ${\mathcal D}_2$ is {\em another dual arc in $A_2(K^{(4)})$   with vector dimensions $(10, 4, 1),$ but it has more members than ${\mathcal D}_1^*$.}

\Example \rm Let $V$ be any vector space over $K$ of dimension $n$.  The exterior algebra 
$ \bigwedge V = K\oplus V\oplus (V\wedge V) \oplus  (V\wedge V\wedge V) \oplus \cdots 
$
is a graded algebra that can replace the polynomial ring in the construction in Section~\ref{A construction}.  
While we obtain   a set of subspaces of a graded component of this algebra with good pairwise intersections, triple intersections 
show that it is no longer a generalized dual arc.

Suppose $\dim V = 5$.  Then
${\mathcal D}_3:=   \{ V\wedge\langle v\rangle \mid 0\neq v\in  V \}  
$ 
 is a set of $1+q+q^2+q^3+q^4$\ \  $4$-spaces of the $10$-space $V\wedge V$.  The intersection of any two members of ${\mathcal D}_3$ is a $1$-space.  
Any $1$-space that is the intersection of two members of ${\mathcal D}_3$ in fact lies in $1+q$ members of ${\mathcal D}_3$, so this is {\em not a dual arc.} 

\medskip
Besides illustrating a  use both of duals and exterior algebra, 
these last three examples possess a numerology that raises a question.  
In view of the fact that the 10-spaces
$A_3(K^{3})^*$, $A_2(K^{4})$ and $K^{5}\wedge K^{5}$
 are isomorphic, is there a relationship among the structures
 ${\mathcal D}_1^*, $ $ {\mathcal D}_2$  and~$ {\mathcal D}_3$?%
 
 \medskip
 {\noindent \em Ackowlegement}:  We are grateful to 
  A. Maschietti and  A. Polishchuk for helpful comments
  concerning Section~\ref{Veronesean points}.

\end{document}